\documentclass[12pt]{amsart}

\usepackage{amssymb}

\usepackage[english,francais]{babel}

\newtheorem{e-proposition}[theorem]{Proposition}

\newtheorem{e-definition}[theorem]{Definition\rm}

\begin{document}

\selectlanguage{english}

\begin{center}

{\huge Lie geometry of flat fronts in hyperbolic space}

\vspace{0.2in}

\selectlanguage{francais}

{\huge La g\'eom\'etrie de Lie des fronts plats dans l'\'espace hyperbolique}

\selectlanguage{english}

\vspace{0.5in}

{\large Francis E Burstall, Udo~Hertrich-Jeromin, Wayne Rossman}

\vspace{0.5in}

\end{center}

\begin{quote}
{\bf Abstract:} We propose a Lie geometric point of view on flat fronts 
in hyperbolic space as special $\Omega$-surfaces and discuss the Lie 
geometric deformation of flat fronts.
\end{quote}

\vspace{0.2in}

\selectlanguage{francais}
\begin{quote}
{\bf R\'esum\'e:} 
Nous proposons un point de vue de Lie g\'eometrie sur les fronts plats
dans l'\'espace hyperbolique comme des surfaces $\Omega$ sp\'eciales.
Nous discutons ensuite la d\'eformation Lie g\'eometrique des fronts
plats.
\end{quote}

\selectlanguage{english}
\def\R{\mathbb{R}}
\def\span{{\rm span}}
\def\proclaim #1.#2\par{\space{\sl#2}\par}

\section{Introduction}
We describe (parallel families of) flat fronts in hyperbolic space
in the realm of Lie sphere geometry:
they turn out to be those $\Omega$-surfaces\footnote{$\Omega$-surfaces
were introduced by Demoulin in \cite{de11a}.} whose enveloped
isothermic sphere congruences each touch a fixed sphere.
This characterization is closely related to the fact that the two
hyperbolic Gauss maps of a flat front are holomorphic,
see \cite{gmm00}.
The fact that ``fronts'' are defined by having a regular Legendre
lift as well as the fact that flat fronts appear in parallel families
both suggest that this is a natural viewpoint.

Since Demoulin's $\Omega$-surfaces can be characterized as envelopes of
isothermic sphere congruences --- such sphere congruences always appearing
in pairs that separate the curvature spheres harmonically --- the rich
theory of isothermic transformations becomes available.
As a first example, we discuss the Calapso transformations of the
isothermic sphere congruences:
these induce the Lie geometric deformation of the surface as a
non-rigid surface in Lie geometry,
and it preserves the condition to project to a flat front in
hyperbolic space --- to see this we employ a pair of linear
conserved quantities that the two fixed spheres give rise to.

Besides the implications for the theory of smooth flat fronts,
our approach also leads to a natural integrable discretization
of flat fronts in hyperbolic space,
cf \cite{hrsy09}.
Thus, we expect this exposition to be the foundation for a wealth
of further research.

{\it Acknowledgements.\/}
We would like to thank M Kokubu, M Umehara and K Yamada for pleasant
and fruitful discussions about the subject.
We also gratefully acknowledge support from the Japan Society for the
Promotion of Science through the second author's fellowship grant
L-08515.

\section{Flat fronts in hyperbolic space as $\Omega$-surfaces}
Let $\mathfrak{f}:M^2\to H^3=\{y\in\R^{3,1}\,|\,|y|^2=-1,y_0>0\}$ be
a flat front in hyperbolic space with unit normal field
$\mathfrak{t}:M^2\to S^{2,1}$.
Away from umbilics and singularities, curvature line coordinates
$(u,v)$ can be introduced,
$\textstyle
   0 = \mathfrak{t}_u + \tanh\varphi\,\mathfrak{f}_u
     = \mathfrak{t}_v + \coth\varphi\,\mathfrak{f}_v
$
with a suitable function $\varphi$,
and by the Codazzi equations $(u,v)$ can be chosen so that
the induced metric takes the form
$\textstyle
   ds^2 = \cosh^2\varphi\,du^2 + \sinh^2\varphi\,dv^2
$
and the Gauss equation becomes harmonicity of $\varphi$.
Now
$$\textstyle
     ({\sqrt{E}\over\sqrt{G}}{\kappa_{1u}\over\kappa_1-\kappa_2})_v
   + ({\sqrt{G}\over\sqrt{E}}{\kappa_{2v}\over\kappa_1-\kappa_2})_u
   = -\varphi_{uv} + \varphi_{vu}
   = 0,
$$
characterizing $\mathfrak{f}$ as an $\Omega$-surface,
see \cite{de11a}.

In order to see this geometrically, consider the orthogonal
decomposition $\R^{4,2}=\R^{1,1}\oplus\R^{3,1}$ and
fix an orthonormal basis $(\mathfrak{p},\mathfrak{q})$ of $\R^{1,1}$,
where $\mathfrak{p}$ defines a point sphere complex, $|\mathfrak{p}|^2=-1$.
Then
$$\textstyle
   (u,v)\mapsto f(u,v) :
   = \span\{\mathfrak{q}+\mathfrak{f}(u,v),\mathfrak{p}+\mathfrak{t}(u,v)\}
$$
defines the Legendre lift of $\mathfrak{f}$ with curvature spheres
$$
   s_1 = \cosh\varphi\,(\mathfrak{p}+\mathfrak{t})
       + \sinh\varphi\,(\mathfrak{q}+\mathfrak{f})
$$
and
$$
   s_2 = \sinh\varphi\,(\mathfrak{p}+\mathfrak{t})
       + \cosh\varphi\,(\mathfrak{q}+\mathfrak{f}),
$$
where the normalizations of the $s_i$ have been chosen so that
$s_{1u}=\varphi_us_2$ and $s_{2v}=\varphi_vs_1$.
Consequently, $s^\pm:=s_1\pm s_2$ are two enveloped isothermic sphere
congruences for $f$ since
$$\textstyle
   s^\pm_{uv} = (\pm\varphi_{uv}+\varphi_u\varphi_v)\,s^\pm,
$$
again characterizing the surface as an $\Omega$-surface,
see \cite{de11a} and \cite{de11b}.
Note that the spheres $s^\pm(u,v)$ separate the curvature spheres
harmonically on the projective line (contact element) given by $f(u,v)$.

As the Legendre map of a flat front is required to be regular
this second analysis and, in particular, the isothermic sphere
congruences $s^\pm$ extend through the singularities of $\mathfrak{f}$.

To characterize those $\Omega$-surfaces that project to flat fronts
in hyperbolic space, first note that the second envelopes of the two
isothermic sphere congruences $s^\pm$ are the hyperbolic Gauss maps
of the flat front:
$\textstyle
   s^\pm \perp \mathfrak{p} \pm \mathfrak{q} =: q^\pm,
$
that is, $s^\pm$ has oriented contact with the infinity sphere
of hyperbolic space $q^\pm$, equipped with opposite orientations.

Thus, suppose $f=\span\{s^+,s^-\}$ is an $\Omega$-surface, given in
terms of a pair of (isothermic) sphere congruences $s^\pm$ that separate
the curvature spheres harmonically.
Assume that each of the sphere congruences $s^\pm$ envelops a fixed
sphere $q^\pm$ so that $q^+$ and $q^-$ do not span a contact element,
that is, $(q^+,q^-)\neq0$; wlog., $(q^+,q^-)=-2$.
Next we fix a point sphere complex
$\textstyle
   \mathfrak{p} := {1\over2}(q^++q^-)
$
and assume that both of the spheres $s^\pm$ never become
point spheres, $(s^\pm,\mathfrak{p})\neq0$.
Then
$$
   \mathfrak{f} :
      = -({q^+\over2~}+{s^+\over(s^+,q^-)})
      +  ({q^-\over2~}+{s^-\over(s^-,q^+)})
$$
and
$$
   \mathfrak{t} :
      = -({q^+\over2~}+{s^+\over(s^+,q^-)})
      -  ({q^-\over2~}+{s^-\over(s^-,q^+)})
$$
take values in $\R^{3,1}=\{q^+,q^-\}^\perp$;
moreover, $|\mathfrak{f}|^2=-1$ so that $\mathfrak{f}$ maps into one of the
hyperbolic spaces with infinity sphere $q^\pm$.
Using the contact condition $(ds^+,s^-)=0$ it is straightforward
to show that $\mathfrak{t}$ is a unit normal field of $\mathfrak{f}$.
Finally, to see that $\mathfrak{f}$ is a flat front, we wheel out the
assumption that $s^\pm$ separate the curvature spheres
$\textstyle
   s_i = -(1+\kappa_i){s^+\over(s^+,q^-)}
       -  (1-\kappa_i){s^-\over(s^-,q^+)}
$
harmonically, that is, the cross ratio
${\kappa_1-1\over\kappa_1+1}{\kappa_2+1\over\kappa_2-1}=-1$
where $\mathfrak{f}$ immerses, implying $\kappa_1\kappa_2=1$.

Note that a different normalization $\tilde q^\pm=e^{\pm t}q^\pm$ of
the $q^\pm$ leads to a parallel flat front.

Thus we have proved\footnote{%
Observe that we did not use that the harmonically separating sphere
congruences $s^\pm$ are isothermic:
indeed an immersed sphere congruence that touches a fixed sphere is
automatically Ribaucour --- hence a Legendre map with two enveloped
sphere congruences that separate the curvature spheres harmonically
and each touches a fixed sphere is automatically an $\Omega$-surface,
cf \cite[\S85]{bl29}.
}:
\proclaim Thm.
Flat fronts in hyperbolic space are those $\Omega$-surfaces with a pair
of isothermic sphere congruences that each touch a fixed sphere, where
the fixed spheres do not span a contact element.
The two fixed spheres yield the two orientations of the infinity sphere
of the hyperbolic ambient space of the flat front;
this determines the point sphere complex, hence the flat front, up to
parallel transformation.

\section{Deformation of flat fronts in hyperbolic space}
Demoulin \cite{de11a} introduced $\Omega$-surfaces as the
Lie geometric analogue of the $R$-surfaces of projective geometry.
It is therefore not too surprising that $\Omega$-surfaces are the
generic\footnote{The other, non-generic class of deformable surfaces,
corresponding to the $R_0$-surfaces, being those where one of the
curvature sphere congruences becomes isothermic in the sense that
it has a Moutard lift, e.g., channel surfaces.}
deformable surfaces of Lie geometry,
see \cite[\S85]{bl29} or \cite{muni06}.
Indeed, each of the isothermic sphere congruences $s^\pm$ enveloped
by an $\Omega$-surface $f$ comes with its
Calapso transformations
$T^\pm(\lambda):M^2\to O(4,2)$,
where $dT^\pm(\lambda)=T^\pm(\lambda)\,\lambda\tau^\pm$
with $\tau^\pm=s^\pm\wedge\star ds^\pm$
for Moutard lifts $s^\pm$ of the isothermic sphere congruences,
$s^\pm\wedge s^\pm_{uv}=0$, and the Hodge-$\star$ operator
of Fubini's quadratic form, $\star du=du$ and $\star dv=-dv$.
Aligning the Moutard lifts to reflect across the Lie cyclides\footnote{%
This is always possible; we have seen it above in the flat front case
that we will discuss.},
$s^\pm=s_1\pm s_2$ for suitable lifts of the curvature spheres,
$$\textstyle
     \tau^+ + {1\over2}d(s^+\wedge s^-)
   = \tau^- + {1\over2}d(s^-\wedge s^+)
   =: \tau
$$
since $s_1\wedge s_{2v}=s_2\wedge s_{1u}=0$,
which leads to
$$\textstyle
   T^+(\lambda)(1+{\lambda\over2}s^+\wedge s^-)
   = T^-(\lambda)(1+{\lambda\over2}s^-\wedge s^+)
   = T(\lambda),
$$
where $T(\lambda)$ yields the Lie-geometric deformation
of $f=\span\{s^+,s^-\}$ via $f(\lambda)=T(\lambda)f$,
see \cite{muni06}.
Observe that, as $f$ is a congruence of null $2$-planes,
$T$ and $T^\pm$ act the same on $s^+$ and $s^-$, hence on $f$:
\proclaim Lemma.
The Calapso transformations $T^\pm$ of the enveloped isothermic sphere
congruences $s^\pm$ of an $\Omega$-surface $f$ yield the Lie geometric
deformation $\lambda\to T(\lambda)f=T^\pm(\lambda)f$ of $f$ as a
Lie deformable surface.

In the flat front case $s^\pm=e^{\pm\varphi}(\mathfrak{p\pm q+t\pm f})$ and
$ds^\pm=\pm d\varphi\,s^\pm-e^{\mp\varphi}\star ds^\mp$ so that
$$\textstyle
   \tau = -{1\over2}\{s^+\wedge ds^-+s^-\wedge ds^+\}
      - d\varphi\,s^+\wedge s^-
      = -(\mathfrak{p+t})\wedge d\mathfrak{t}
      + (\mathfrak{q+f})\wedge d\mathfrak{f}
$$
and the two fixed sphere congruences $q^\pm$ give rise to two linear
conserved quantities
$$\textstyle
   p^\pm(\lambda)
   = (1+{\lambda\over2}s^\pm\wedge s^\mp)^{-1}q^\pm
   = q^\pm+{\lambda\over2}(q^\pm,s^\mp)s^\pm
$$
$$
   = (1-\lambda)(\mathfrak{p\pm q}) - \lambda(\mathfrak{t\pm f})
$$
for the Lie geometric deformation $T(\lambda)$ since
$d(Tp^\pm)(\lambda)=T(\lambda)(d+\lambda\tau)p^\pm(\lambda)=0$.
Hence, the deformed surfaces are flat fronts in hyperbolic space as long
as $(p^+,p^-)(\lambda)=-2(1-2\lambda)\neq0$ with $(Tp^\pm)(\lambda)$ as
the infinity sphere of the hyperbolic ambient space with its two
orientations:
\proclaim Thm.
The Lie geometric deformation of a flat front in hyperbolic space
yields a $1$-parameter family of flat fronts in hyperbolic space.

Normalizing the conserved quantities,
$p^\pm(\lambda)\to{p^\pm(\lambda)\over\sqrt{1-2\lambda}}$,
and following the construction in the previous section we obtain
$$\textstyle
   \mathfrak{f}(\lambda)
   = T(\lambda){h^+(\lambda)-h^-(\lambda)\over2\sqrt{1-2\lambda}}
      \quad{\rm and}\quad
   \mathfrak{t}(\lambda)
   = T(\lambda){h^+(\lambda)+h^-(\lambda)\over2\sqrt{1-2\lambda}}
$$
with $h^\pm(\lambda)
   = -\lambda(\mathfrak{p\pm q})+(1-\lambda)(\mathfrak{t\pm f})$.
Choosing constants of integration so that $Tp^\pm\parallel q^\pm$ for
every $\lambda\neq{1\over2}$, the deformation is confined to
$\R^{3,1}=\{q^+,q^-\}^\perp$.
With
$\mathfrak{e}_1:=\pm e^{\pm\varphi}(\mathfrak{t\pm f})_u$
and
$\mathfrak{e}_2:=e^{\pm\varphi}(\mathfrak{t\pm f})_v$
we obtain frames
$F(\lambda):=T(\lambda)(\mathfrak{e}_1,\mathfrak{e_2},
   {h^+(\lambda)\over\sqrt{1-2\lambda}},
   {h^-(\lambda)\over\sqrt{1-2\lambda}})$
for the family of point-pair maps into the conformal $2$-sphere
given by the two hyperbolic Gauss maps of the $\mathfrak{f}(\lambda)$.
As (degenerate) Darboux pairs these are curved flats in the space
of point pairs, see \cite[\S5.5.20 or Sect 8.7]{imdg}.
Now
$$
  (d+\lambda\tau)h^\pm(\lambda) = 
  \sqrt{1-2\lambda}\,e^{\mp\varphi}(\pm\mathfrak{e}_1du+\mathfrak{e}_2dv)
$$
$$
   (d+\lambda\tau)\mathfrak{e}_1 
      = (\varphi_vdu-\varphi_udv)\,\mathfrak{e}_2
      + {\sqrt{1-2\lambda}\over2}
         (e^\varphi h^+(\lambda)-e^{-\varphi}h^-(\lambda))\,du
$$
$$
   (d+\lambda\tau)\mathfrak{e}_2 
      = (-\varphi_vdu+\varphi_udv)\,\mathfrak{e}_1
      - {\sqrt{1-2\lambda}\over2}
         (e^\varphi h^+(\lambda)+e^{-\varphi}h^-(\lambda))\,dv
$$
showing that $T(\lambda)({h^+(\lambda)\over\sqrt{1-2\lambda}},
   {h^-(\lambda)\over\sqrt{1-2\lambda}})$
is the $1$-parameter family of Darboux pairs obtained from the
curved flat associated family \cite[\S3.3.3]{imdg}:
\proclaim Thm.
The Lie geometric deformation of a flat front in hyperbolic space
with parameter $\lambda$ yields the curved flat associated family
with parameter\footnote{Note that an imaginary parameter still yields
a real geometry, as the semi-Riemannian symmetric space of point-pairs
in a conformal $n$-sphere is self dual; a change of sign of the curved
flat parameter does not change the geometry of the curved flat,
see \cite[\S5.5.19]{imdg}.}
$\sqrt{1-2\lambda}$ of the Darboux pair of its hyperbolic Gauss
maps in the symmetric space of point-pairs in $S^2$.

Note that a Darboux pair in $S^2$ gives rise to a unique parallel
family of flat fronts in hyperbolic space by determining the orthogonal
surfaces of the cyclic system given by the circles intersecting the
infinity sphere of hyperbolic space orthogonally in the points of the
hyperbolic Gauss maps\footnote{This, in fact, shows that flat fronts
are Guichard surfaces, which were shown to be $\Omega$-surfaces
in \cite{de11a}.}, see \cite{krsuy04}.
Thus the family of Darboux pairs is sufficient to determine the
corresponding Lie geometric deformation of (parallel families of)
flat fronts.

In the case of the peach front, see \cite{kruy05},
the curved flat system can explicitly be integrated:
when $\sqrt{1-2\lambda}$ is real, the $\mathfrak{f}(\lambda)$ become
the snowman type flat fronts, see \cite{krsuy04},
or a new type of flat fronts whereas, when the parameter is imaginary,
we obtain the hourglass type flat fronts, see \cite{krsuy04}.

\end{document}